\documentstyle[amscd,amssymb,verbatim,11pt]{amsart}

\theoremstyle{plain}
\newtheorem{Prop}{Proposition}[section]
\newtheorem{Thm}[Prop]{Theorem}
\newtheorem{Cor}[Prop]{Corollary}

\newtheorem{Lem}[Prop]{Lemma}

\theoremstyle{definition}
\newtheorem{Def}[Prop]{Definition}

\theoremstyle{remark}
\newtheorem{Rem}[Prop]{Remark}
\newtheorem{Problem}[Prop]{\bf Problem}

\def\dim{\mathop{\roman{dim}}}
\def\int{\mathop{\roman{int}}}

\def\1{^{-1}}

\def\dim{\text{dim}}

\def\int{\text{Int}}
\numberwithin{equation}{section}
\begin{document}
\title[
Coarse dimensions and partitions of unity
]%
   {Coarse dimensions and partitions of unity
}

\author{N.~Brodskiy}
\address{University of Tennessee, Knoxville, TN 37996, USA}
\email{brodskiy@@math.utk.edu}

\author{J.~Dydak}
\address{University of Tennessee, Knoxville, TN 37996, USA}
\email{dydak@@math.utk.edu}

\date{ September 15, 2005}
\keywords{Asymptotic dimension, coarse dimension, coarse category,
Lebesque number}

\subjclass{ Primary: 54F45, 54C55, Secondary: 54E35, 18B30, 54D35, 54D40, 20H15}

\thanks{ 
}

\begin{abstract} Gromov \cite{Gr$_1$} 
 and Dranishnikov  \cite{Dr$_1$} introduced asymptotic and coarse
dimensions of proper metric spaces via quite different ways.
We define coarse and asymptotic dimension of all
metric spaces in a unified manner and we investigate relationships between them
generalizing results of Dranishnikov  \cite{Dr$_1$} and Dranishnikov-Keesling-Uspienskij
 \cite{DKU}. 
\end{abstract}

\maketitle

\medskip
\medskip
\tableofcontents

\section{Introduction}

There are three concepts of dimension associated with variants of the coarse category
of proper metric spaces.
The original one, the asymptotic dimension of Gromov \cite{Gr$_1$}, and dimensions
$asdim^\ast(X)$ and $\dim^c(X)$ introduced by Dranishnikov \cite{Dr$_1$}.
All three dimensions are defined in seemingly different ways:

1. The asymptotic dimension of Gromov (see \cite{Gr$_1$} or Definitions 1-2 in \cite{Dr$_1$}
on p.1103) is the smallest integer $n$
such that for every $M > 0$ there is a uniformly bounded family
${\cal {U}}$ of Lebesque number at least $M$ and multiplicity (or order)
at most $n+1$.

2. The asymptotic dimension $asdim^\ast(X)$ of Dranishnikov
(see Definition 3 in \cite{Dr$_1$}
on p.1104) is the smallest integer $n$
such that for every proper function $f:X\to R_+$
there is a contracting map $\phi:X\to K$ to an $n$-dimensional asymptotic polyhedron
such that for each $M > 0$ there is a compact subset $C$ of $X$
with the property that $\phi^{-1}(B(\phi(x),M))\subset B(x,f(x))$
for all $x\in X\setminus C$.

3. The coarse dimension $\dim^c(X)$ of Dranishnikov
(see Definition 4 in \cite{Dr$_1$}
on p.1105) is the smallest integer $n$
such that $R^{n+1}$ is an absolute extensor of $X$ in the category
of proper asymptotically Lipschitz functions. That dimension coincides
with the dimension of the Higson corona $\nu(X)$ of $X$
(see Theorem 6.6 in \cite{Dr$_1$}
on p.1111).

One of the main motivations behind the research in asymptotic dimension is the result
of Yu (see  \cite{Yu$_1$} and \cite{Yu$_2$}) that the Novikov Conjecture
holds for groups of finite asymptotic dimension.

In this paper we work in the coarse category of all metric spaces
and we devise a unified way of defining five dimensions:
{\it coarse dimension} $\dim^{coa}_{rse}(X)$,
{\it major coarse dimension} $\dim^{COA}_{RSE}(X)$, {\it asymptotic dimension} $asdim(X)$, {\it minor asymptotic dimension} $ad(X)$,
and {\it large scale dimension} $\dim^{large}_{scale}(X)$. 

In case of proper metric spaces, three of them coincide
with the above dimensions. Namely, $\dim^{COA}_{RSE}(X)=dim^\ast(X)$,
$\dim^{coa}_{rse}(X)=\dim^c(X)$, and $asdim(X)$ coincides with Gromov's asymptotic dimension.
The fourth one, the minor asymptotic dimension, is a variant
of Gromov's dimension. The large scale dimension is always equal to the coarse dimension and the reason we are introducing it is to simplify proofs of the relations between the three basic dimensions which we do
in a much simpler way than as described in Dranishnikov's paper \cite{Dr$_1$}.
The main relations between dimensions are as follows:
\begin{enumerate}
\item There are two strands of inequalities: 
$asdim(X)\ge \dim^{COA}_{RSE}(X)\ge \dim^{coa}_{rse}(X)$
and $asdim(X)\ge ad(X)\ge \dim^{coa}_{rse}(X)$,
\item In each strand (for unbounded spaces $X$), finiteness
of a larger dimension implies its equality with all smaller dimensions in the strand.
\end{enumerate}

We do not know of any unbounded space $X$ such that a larger dimension in a strand
is infinite and a smaller dimension is finite.

Our fundamental concept is that of a coarse family and we follow the well-established route
of defining the covering dimension by refining covers with covers of a prescribed multiplicity.
In classical dimension theory one deals with two cases: finite covers and infinite covers.
There, for paracompact spaces, the two concepts coincide. In the case of coarse covers we get two
concepts of coarse dimension whose equality remains unresolved.
\par
A finite family ${\cal {U}}$ of subsets of $X$ is coarse if and only if there
is a slowly oscillating partition of unity $f$ on $X\setminus B$
for some bounded subset $B$ of $X$ whose carriers $Carr(f)$ refine ${\cal {U}}$.
That explains why, in the case of a proper metric space $X$,
its coarse dimension equals the covering dimension of
the Higson corona of $X$.
\par
Our basic strategy is to associate natural functions with objects
and declare those objects to be coarse, asymptotic, or large scale
if the function is coarsely proper. A function $f$ is {\it coarsely proper}
 if $f(E_n)\to\infty$
whenever $E_n\to\infty$. Elements $E_n$ related to objects could be points in a metric space,
bounded subsets in a metric space, or covers of a metric space
(in which case divergence to infinity is measured by the size of
the Lebesque number). In \cite{Dr$_1$} (p.1089) coarsely proper functions were defined
as those $f:X\to Y$ such that $f^{-1}(A)$ is bounded whenever $A$ is bounded in $Y$.
Notice that our definition generalizes the one from \cite{Dr$_1$}.

\par The authors are grateful to Jose Higes for helpful comments.

\section{Preliminaries}

Given a subset $A\ne\emptyset$ of a metric space $X$ the most basic function
is {\it the distance function}
$d_A:X\to R_+$: $d_A(x)=dist(x,A)$.

\begin{Def} \label{BallsDef}
Given a subset $A$ of a metric space $(X,d_X)$ the ball $B(A,M)$
is defined to be the set $\{x\in X \mid  dist(x,A) < M\}$ if $M > 0$,
 it is defined to be the set 
$\{x\in X \mid  dist(x,X\setminus A) > -M\}$ if $M < 0$,
and it is simply $A$ if $M=0$.
\end{Def}

The distance function leads to the first concept of divergence to infinity:
$x_n\to\infty$ if $d_X(x_n,x_0)\to\infty$ for some (and hence for all) $x_0\in X$.
However, $dist(x,A)$ is a function of two arguments and we can use the second
one to define divergence to infinity for bounded subsets of $X$.
Here is a more general concept.

\begin{Def} \label{ProperFamilyDef}
A family ${\cal {U}}$ of bounded subsets of $X$ is called {\it coarsely proper}
if the function $U\to d_U(x_0)$ is coarsely proper for some (and hence for all) $x_0\in X$.
Here ${\cal {U}}$ is considered
as a subspace of all bounded subsets of $X$ with the Hausdorff metric.
\end{Def}

Notice that a sequence $\{A_n\}$ of bounded subsets of $X$
containing points $x_n\in A_n$ so that $x_n\to\infty$ is coarsely proper
if and only if every bounded subset of $X$ intersects at most finitely many
elements of the sequence. In that case we write $A_n\to\infty$ and that
form of divergence to infinity is of most interest to us.

\begin{Lem} \label{SelectionFunAreProper}
If ${\cal {U}}$ is a coarsely proper cover of $X$, then every selection function
$\phi:X\to {\cal {U}}$ (that means $x\in \phi(x)$) is coarsely proper.
\end{Lem}
\begin{pf} Suppose $x_n\to\infty$ and $x_n\in U_n\in {\cal {U}}$.
Clearly, $U_n\to\infty$ in the Hausdorff metric.
Pick $x_0\in X$. Since $d_{U_n}(x_0)\to\infty$, every bounded subset of $X$ 
intersects at most finitely many
elements of the sequence $\{U_n\}$ and any selection function $\phi$ is coarsely proper.
\end{pf}

\begin{Def} \label{LocalLebesqueNumDef}
Given a family ${\cal {U}}$ in $X$, {\it the local Lebesque number}
$L_{{\cal {U}}}(x)\in R_+\cup\infty$ is defined as the supremum
of $dist(x,X\setminus U)$, $U\in{\cal {U}}$. If $U=X$ for some $U\in{\cal {U}}$
it is defined to be infinity.
\end{Def}

Notice that either $L_{{\cal {U}}}\equiv\infty$ at all points or
it is a natural Lipschitz function associated with ${\cal {U}}$.
More precisely $|L_{{\cal {U}}}(x)-L_{{\cal {U}}}(y)|\leq d_X(x,y)$.

\begin{Def} \label{LebesqueNumDef}{\it The Lebesque
number} $L({\cal {U}},A)$ is $\inf\{L_{{\cal {U}}}(x) \mid x\in A\}$.
\end{Def}

\begin{Def} \label{CoarseFamilyDef}
A family of subsets ${\cal {U}}$ of a metric space $X$
is called {\it coarse} if $L_{{\cal {U}}}$ is coarsely proper
(as a function from $X$ to $R\cup\infty$).
\end{Def}

An alternative way to define coarse families is to require
$L({\cal {U}},A)\to\infty$ as $A\to\infty$. Yet another way
is to state that $L({\cal {U}},X\setminus B(x_0,t))\to\infty$
as $t\to\infty$.

\begin{Prop}\label{CharOfFiniteCoarseFam}
\begin{itemize}
\item[1.] A family ${\cal {U}}=\{A\}$ consisting of one subset $A$ of $X$
is coarse if and only if $X\setminus A$ is bounded.

\item[2.] A family ${\cal {U}}=\{X_1,X_2\}$ consisting of two subsets
of $X$ is coarse if and only if $d_X$ 
restricted to $(X\setminus X_1)\times (X\setminus X_2)$ is coarsely proper.
\item[3.] A family ${\cal {U}}=\{X_1,X_2,\ldots,X_n\}$ consisting of finitely many subsets
of $X$ is coarse if and only if the function $d_{{\cal {U}}}(x):=\sum\limits_{i=1}^n dist(x,X\setminus X_i)$
is coarsely proper.
\end{itemize}
\end{Prop}
\begin{pf} 1. If $X\setminus A$ is bounded, then $L_{{\cal {U}}}(x)\ge dist(x,X\setminus A)$ 
and $L_{{\cal {U}}}$ is coarsely proper. If $X\setminus A$ is unbounded, then
$L_{{\cal {U}}}(x)=0$ at all $x\in X\setminus A$ and $L_{{\cal {U}}}$ is not coarsely proper.

2. Suppose ${\cal {U}}=\{X_1,X_2\}$ is coarse 
and $x_n\to\infty$, $y_n\to\infty$,
for some $x_n\in X\setminus X_1$, $y_n\in X\setminus X_2$.
Notice $L_{{\cal {U}}}(x_n)\leq d_X(x_n,y_n)$, so $d_X(x_n,y_n)\to\infty$.

If ${\cal {U}}=\{X_1,X_2\}$ is not coarse, then there is a sequence
$z_n\to\infty$ with $L_{{\cal {U}}}(z_n)$ bounded by $M$.
We can produce $x_n\in X\setminus X_1$ and
$y_n\in X\setminus X_2$ so that $d_X(z_n,x_n) < M+1$ and
$d_X(z_n,y_n) < M+1$ for all $n$. Now, $d_X(x_n,y_n) < 2M+2$,
a contradiction.

3. Notice $d_{{\cal {U}}}(x) \ge L_{{\cal {U}}}(x)$
and $m\cdot L_{{\cal {U}}}(x)\ge d_{{\cal {U}}}(x)$.
\end{pf}

\begin{Def} \label{LebesqueNumTransferDef}
Given a function $f:X\to Y$ of metric spaces, its {\it Lebesque number transfer}
$L^f:R_+\to R_+\cup\infty$ is the supremum of all
functions $\alpha :R_+\to R_+\cup\infty$ such that
$L({\cal {U}},Y) \ge t$ implies
$L(f^{-1}({\cal {U}}),X)\ge \alpha (t)$ for all families ${\cal {U}}$
of subsets of $Y$.
\end{Def}

\begin{Def}\label{CoarseFunctionDef}
A function $f:X\to Y$ of metric spaces is {\it coarse}
if the Lebesque number transfer $L^f$ is coarsely proper.
\end{Def}

An alternative definition of coarse functions
is to require the function ${\cal {U}}\to L(f^{-1}({\cal {U}}),X)$
to be coarsely proper on the set of covers of $Y$.

Let us show that our definition of coarse functions coincides with that of Roe
\cite{R$_2$}.

\begin{Prop}\label{CoarseAndRoeCoarse}
A function $f:X\to Y$ is coarse if and only if
for every $R > 0$ there is $M > 0$ such that $d_X(x,y)\leq R$
implies $d_Y(f(x),f(y))\leq M$ for all $x,y\in X$.
\end{Prop}
\begin{pf} Notice that if $M > 0$ and $N > 0$
are numbers such that $d_X(x,y) < M$ implies $d_Y(f(x),f(y)) < N$,
then $L^f(N)\ge M$. Therefore $f$ being coarse in the sense of Roe implies
$L^f$ being coarsely proper.

Conversely, if $L^f(N)\ge M$, then consider the cover ${\cal {U}}=\{B(z,N)\}_{z\in Y}$
whose Lebesque number is clearly at least $N$.
If $d_X(x,y) < M$, then there is $z$ so that $x,y\in f^{-1}(B(z,N))$.
Hence $d_Y(f(x),f(y)) < 2\cdot N$ and $f$ is coarse.
\end{pf}

Dranishnikov \cite{Dr$_1$} (p.1088) defined {\it asymptotically Lipschitz functions}
$f:X\to Y$ as those for which there are constants $M > 0$ and $A$ such that
$d_Y(f(x),f(y))\leq M\cdot d_X(x,y)+A$ for all $x,y\in X$.
Let us relate this concept to the Lebesque number transfer.
\begin{Prop}\label{AsyLipschitzChar}
A function $f:X\to Y$ is asymptotically Lipschitz if and only if
there is a linear function $t\to m\cdot t+b$ so that
$m > 0$ and $L^f(t)\ge m\cdot t+b$ for all $t$.
\end{Prop}
\begin{pf} Suppose there are constants $M > 0$ and $A$ such that
$d_Y(f(x),f(y))\leq M\cdot d_X(x,y)+A$ for all $x,y\in X$.
Given a cover ${\cal {U}}$ of $Y$ with $L({\cal {U}},Y)\ge t$
and given $x\in X$, the ball $B(x,(t-A-\delta )/M)$
is mapped by $f$ into the ball $B(f(x),t-\delta )$ which is contained
in an element of ${\cal {U}}$ for all $\delta >0$.
That shows the Lebesque number of $f^{-1}({\cal {U}})$ to be at least
$(t-A)/M)$.
Conversely, if $L^f(t)\ge m\cdot t+b$ for all $t$ and $m > 0$,
then we claim $d_Y(f(x),f(y)) < 2\cdot d_X(x,y)/m+2(1-b)/m$.
Indeed, put $d_X(x,y)=s$ and consider the cover ${\cal {U}}=\{B(z,(s+1-b)/m)\}_{z\in Y}$
whose Lebesque number is clearly at least $(s+1-b)/m$.
There is $z$ so that $x,y\in f^{-1}(B(z,(s+1-b)/m))$.
Hence $d_Y(f(x),f(y)) < 2\cdot (s+1-b)/m$ and $f$ is asymptotically Lipschitz.
\end{pf}

\begin{Prop}\label{ProperCoarseAndPtInv}
Given a function $f:X\to Y$ of metric spaces
the following conditions are equivalent:
\begin{itemize}
\item[1.] $f$ sends bounded subsets of $X$
to bounded subsets of $Y$ and
$f^{-1}({\cal {U}})$ is coarse for every coarse family ${\cal {U}}$ in $Y$.

\item[2.] $f$ is coarse and coarsely proper.
\end{itemize}
\end{Prop}
\begin{pf} 1$\implies$2. Given a bounded subset $A$ of $Y$
the family $\{Y\setminus A\}$ is coarse (see \ref{CharOfFiniteCoarseFam}).
Since $\{f^{-1}(Y\setminus A)\}$ is coarse and
$f^{-1}(Y\setminus A)=X\setminus f^{-1}(A)$, $f^{-1}(A)$ must be bounded
and $f$ is coarsely proper.

If $f$ is not coarse, we find sequences $x_n,y_n\in X$ so that
$d_Y(f(x_n),f(y_n)) > n$ for each $n$ but $d_X(x_n,y_n) < M$ for all $n$.
Since $f$ sends bounded subsets of $X$
to bounded subsets of $Y$, we may assume $x_n\to\infty$, hence $y_n\to\infty$.
Put $A=\{x_n\}$ and $B=\{y_n\}$.
Using \ref{CharOfFiniteCoarseFam} we see that ${\cal {U}}=\{Y\setminus f(A), Y\setminus f(B)\}$
is a coarse family in $Y$. Since $f^{-1}({\cal {U}})$ is coarse,
the family ${\cal {V}}=\{X\setminus A, X\setminus B\}$, to which ${\cal {U}}$
is a shrinking, is coarse as well. That however contradicts \ref{CharOfFiniteCoarseFam}.

2$\implies$1. Obviously, coarse functions $f:X\to Y$ send bounded subsets of $X$
to bounded subsets of $Y$.
Put ${\cal {V}}=f^{-1}({\cal {U}})$ for some coarse family
${\cal {U}}$ in $Y$. To find points $x\in X$ such that $L_{{\cal {V}}}(x) > t$
we find $s > 0$ so that $L^f(s) > t$ and we find $u > 0$ such that
$L_{{\cal {U}}}(y) > s$ for $y\in Y\setminus B(y_0,u)$.
Put ${\cal {W}}={\cal {U}}\cup \{B(y_0,u+s)\}$. Note $L({\cal {W}},Y) > s$.
Since $L(f^{-1}({\cal {W}}),X) > t$, points $x$ lying outside of the bounded set
$f^{-1}(B(y_0,u+s))$ satisfy $L_{{\cal {V}}}(x) > t$.
\end{pf}

In the end of this section let us demonstrate the usefulness
of the concept of a coarse family by rewording notions from
\cite{DZ}.

In section 5.2 of \cite{DZ} the concept
of {\it asymptotic neighborhood} $W$ of a subset $A$ of $X$
is introduced by requiring $\lim\limits_{r\to\infty}dist(A\setminus B(x_0,r),X\setminus W)=\infty$
for some (and hence for all) $x_0\in X$.

\begin{Prop}\label{AsyNbhdViaCoarse}
$W$ is an asymptotic neighborhood of $A$ if and only if
the pair $\{X\setminus A,W\}$ is coarse. 
\end{Prop}
\begin{pf} According to part 2 of \ref{CharOfFiniteCoarseFam} the pair
$\{X\setminus A,W\}$ is coarse if and only if
$d_X$ restricted to $A\times (X\setminus W)$ is coarsely proper. That can be
easily seen as equivalent to $\lim\limits_{r\to\infty}dist(A\setminus B(x_0,r),X\setminus W)=\infty$
for some (and hence for all) $x_0\in X$.
\end{pf}

In section 5.2 of \cite{DZ} (see also \cite{Dr$_2$}) the concept
of {\it asymptotically disjoint subsets} $A$ and $B$ of $X$
is introduced by requiring $\lim\limits_{r\to\infty}dist(A\setminus B(x_0,r),B\setminus B(x_0,r))=\infty$
for some (and hence for all) $x_0\in X$.

\begin{Prop}\label{AsyDisjViaCoarse}
$A$ and $B$ are asymptotically disjoint if and only if
the pair $\{X\setminus A,X\setminus B\}$ is coarse. 
\end{Prop}
\begin{pf} Apply part 2 of \ref{CharOfFiniteCoarseFam}.
\end{pf}

Also notice that the concept of an asymptotic separator  
 of \cite{DZ} (see section 5.2)
can be introduced without referring to the Higson corona.

\begin{Def}\label{AsySeparator}
A subset $C$ of $X$ is an {\it asymptotic separator}
between asymptotically disjoint subsets $A$ and $B$
if there are asymptotic neighborhoods $W_A$ of $A$
and $W_B$ of $B$ such that $C=X\setminus (W_A\cup W_B)$
and $W_A\cap W_B=\emptyset$.
\end{Def}

\section{Multiplicity and higher Lebesque numbers}

\begin{Def} \label{MultDef}
Given a family ${\cal {U}}$ of subsets of $X$ we define
{\it the multiplicity function}
$m_{{\cal {U}}}:X\to Z_+\cup\infty$ by setting $m_{{\cal {U}}}(x)$
to be equal to the number of elements of ${\cal {U}}$ containing $x$.
The {\it global multiplicity} $m({\cal {U}},A)$ is the supremum
of $m_{{\cal {U}}}(x)$, $x\in A$.
\end{Def}

By a {\it coarse refinement} ${\cal {V}}$ of a coarse
family ${\cal {U}}$ we mean a coarse family such that
every element $V$ of ${\cal {V}}$ is contained in an
element $U$ of ${\cal {U}}$.
${\cal {V}}$ is called a {\it shrinking} of ${\cal {U}}$ if they are indexed
by the same set $S$ and $V_s\subset U_s$ for all $s\in S$.
If ${\cal {V}}$ is a coarse refinement of ${\cal {U}}$ indexed by a different
set $T$, then one can create a shrinking ${\cal {V}}'$ of ${\cal {U}}$
as follows: find a function $\phi:T\to S$ satisfying $V_t\subset U_{\phi(t)}$
for all $t\in T$. Define $V'_s$ as $\bigcup\{V_t \mid  s=\phi(t)\}$.
Notice that ${\cal {V}}'$ has multiplicity at most that of ${\cal {V}}$
and is a coarse shrinking of ${\cal {U}}$.
\par

Given a family $\phi=\{\phi_s:X\to R_+\}_{s\in S}$
of functions its {\it carrier family} $Carr(\phi)$ is the family
$\{\phi_s^{-1}(0,\infty)\}_{s\in S}$.
The {\it multiplicity} $m(\phi)$ of $\phi$ is defined as the multiplicity
of its carrier family and its {\it Lebesque number} $L(\phi)$
is defined as the Lebesque number of its carrier family.

\begin{Lem}\label{CoversWithInfiniteLebesque}
If ${\cal {U}}=\{U_s\}_{s\in S}$ is family in $X$ such that
$L_{{\cal {U}}}(x_0)=\infty$ for some $x_0\in X$, then it has a coarse
refinement ${\cal {V}}$ of multiplicity at most $2$.
\end{Lem}
\begin{pf} 
Put $V_{n}=\{x\in X \mid  (n-1)^2\leq d(x,x_0) < (n+1)^2\}$ for $n\ge 1$.
\end{pf}

\begin{Lem} \label{OstrandPropertyForCovers}
If ${\cal {U}}=\{U_s\}_{s\in S}$ is a family
in $X$ of multiplicity at most $n+1$,
then it can be refined by ${\cal {V}}=\bigcup\limits_{i=1}^{n+1}{\cal {V}}^i$
such that $L_{{\cal {V}}}(x)\ge L_{{\cal {U}}}(x)/(2n+2)$ for each $x\in X$
and each ${\cal {V}}^i$ consists of disjoint sets.
\end{Lem}
\begin{pf} 
Define $f_s(x)=dist(x,X\setminus V_s)$. For each finite set $T$ of $S$
define $W_T=\{x\in X \mid  \min\{f_t(x)\mid  t\in T\} > \sup\{f_s(x) \mid  s\in S\setminus T\}\}$.
Notice $W_T=\emptyset$ if $T$ contains at least $n+2$ elements.
Also, notice that $W_T\cap W_F=\emptyset$ if both $T$ and $F$ 
are different but contain the same
number of elements.
Let us estimate the Lebesque number of ${\cal {W}}=\{W_T\}_{T\subset S}$.
Given $x\in X$ arrange all non-zero values $f_s(x)$ from the largest
to the smallest. Add $0$ at the end and look at gaps between those values.
The largest number is at least $L_{{\cal {U}}}(x)$, there are at most $n+1$ gaps,
so one of them is at least $L_{{\cal {U}}}(x)/(n+1)$. That implies the ball $B(x,L_{{\cal {U}}}(x)/(2n+2))$
is contained in one $W_T$ ($T$ consists of all $t$ to the left of the gap).
Define ${\cal {V}}_i$ as $\{W_T\}$, all $T$ containing exactly $i$ elements.
\end{pf}

\begin{Lem}\label{ProperRefinementsOfCoarse}
If \ ${\cal {U}}=\{U_s\}_{s\in S}$ is a coarse family in $X$, then it has a coarse
refinement ${\cal {V}}$ that is coarsely proper.
Moreover, if \ ${\cal {U}}$ is of finite multiplicity,
then we may require ${\cal {V}}$ to be of finite multiplicity
as well.
\end{Lem}
\begin{pf}
Let ${\cal {V}}=\{V_{s,m}\}_{(s,m)\in S\times N}$,
where $V_{s,m}=\{x\in U_s \mid  2^m< d(x,x_0)\leq 2^{m+2}\}$.
Notice ${\cal {V}}$ is coarse of multiplicity at most $2\cdot m({\cal {U}})$.
Also, it consists
of bounded sets so that for any sequence $x_k\to\infty$
the conditions $x_k\in V_{s(k),m(k)}$ imply $V_{s(k),m(k)}\to\infty$.
\end{pf}

\begin{Prop}\label{CoarseShrinkingALaParacompact}
If ${\cal {U}}=\{U_s\}_{s\in S}$ is a coarse family in $X$, then it has a coarse
shrinking ${\cal {V}}=\{V_s\}_{s\in S}$ such that for any $M > 0$
there is a bounded subset $A_M$ of $X$ with the property that $B(x,M)\cap V_s\ne \emptyset$
implies $B(x,M)\subset U_s$ provided $x\in X\setminus A_M$.
\end{Prop}
\begin{pf}
Pick $x_0\in X$ and define $f(x)=\min(d(x,x_0)/2,L_{{\cal {U}}}(x)/2)$.
Notice $f$ is a coarsely proper function of Lipschitz constant $1/2$.
For each $x\in X$ pick $s(x)\in S$ so that $B(x,f(x))\subset U_{s(x)}$.
Define $V_s$ as the union of those balls $B(x,f(x)/2)$ so that $s=s(x)$.
It suffices to observe that $B(x,M)\cap V_s\ne\emptyset$ and $M < f(x)/3$
implies $B(x,M)\subset U_s$. Indeed, $B(y,f(y))\subset U_s$ for some
$y\in B(x,M)$. Since $f(x)-f(y) \leq d(x,y)/2 < M/2$,
one has $f(y) > f(x)-M/2> 3M-M/2 > 2M$ and $B(x,M)\subset B(y,f(y))\subset U_s$.
\end{pf}

\begin{Lem}\label{CoarseShrinkingWithProperFun}
If ${\cal {U}}$ is a coarse family in $X$ that is coarsely proper, then there is
a coarsely proper function $f:{\cal {U}}\to R_+$ such that
the family $\{B(U,-f(U))\}_{U\in {\cal {U}}}$ is coarse.
\end{Lem}
\begin{pf}
Define $f(U)=\inf\{L_{{\cal {U}}}(x)/4 \mid  x\in U\}$.
Notice $f$ is a coarsely proper function.
Pick $s(x)\in S$ so that $B(x,L_{{\cal {U}}}(x)/2)\subset U_{s(x)}$.
$f(U_{s(x)}) \leq L_{{\cal {U}}}(x)/4$ which implies
$B(x,L_{{\cal {U}}}(x)/4)\subset  B(U_{s(x)},-f(U_{s(x)}))$.
Thus $\{B(U,-f(U))\}_{U\in {\cal {U}}}$ is coarse.
\end{pf}

In the large scale geometry one should think of bounded subsets of $X$
as points. Here is a generalization of the Lebesque number.

\begin{Def}\label{HigherLebesqueNumDef} Let $n\ge 0$. Suppose ${\cal {U}}$ is a family in $X$ and $A$ is
a bounded subset of $X$. The {\it $n$-th Lebesque number} $L^{n}({\cal {U}},A)$
 is the supremum of $t\in [0,\infty]$ such that
${\cal {U}} |_  A$ has a refinement of multiplicity at most $n+1$ and Lebesque
number at least $t$.
\end{Def}

Notice such supremum exists as 
the cover of $A$ consisting of points is of Lebesque number $0$ and
multiplicity $1$.

Observe that $L^{n}({\cal {U}},A)$, $n\ge 0$, form an increasing sequence
of numbers bounded by $L({\cal {U}},A)$. If ${\cal {U}}|_ A$ is of finite order, then
they eventually stabilize and are equal to $L({\cal {U}}|_ A,A)$.

Let us point out that Sperner's Lemma can be used to estimate higher Lebesque
numbers as follows: Consider a 2-simplex $\Delta$ with vertices labeled
$0$, $1$, and $2$. Let ${\cal {U}}$ be the cover of $\Delta$ by stars $U_i$, $i=0,1,2$, of its
vertices. Consider a subdivision $L$ of $\Delta$ with mesh $M$ (in this case it coincides
with the longest edge in the subdivision). Let $X=A$ be the set of vertices of $L$.
Suppose ${\cal {V}}=\{V_0,V_1,V_2\}$ is a shrinking of ${\cal {U}}|_ A$. Obviously, there is a shrinking
of multiplicity $1$. However, if we request ${\cal {V}}$ to be of large Lebesque
number, we run into problems. Namely, $L^{1}({\cal {U}},A) \leq M$.
Indeed, if $L({\cal {V}}) > M$, we assign to each vertex $v$ of $L$ number $i$
such that $v\in V_i$. We are in the situation of the classical Sperner's Lemma:
vertices on the edges of $\Delta$ must be labeled with a number of one
of the vertices of that edge. Therefore one has a simplex in $L$ whose
vertices were assigned all three numbers $0,1,2$. Since $L({\cal {V}}) > M$,
the three vertices belong to $V_0\cap V_1\cap V_2$ and multiplicity
of ${\cal {V}}$ is $3$. Thus $L^{1}({\cal {V}},A) \leq M$.
\par We will use the observation above in the case of $M$-scale
connected spaces.

\begin{Def}\label{MScaleConnectedDef} 
Suppose $M$ is a positive number. A metric space $X$ is called
{\it $M$-scale connected} if for every two points $x,y\in X$
there is a chain of points $x=x_1, x_2,\ldots, x_k=y$
such that $d_X(x_i,x_{i+1})< M$ for all $i < k$.
\end{Def}

Here is an application of Sperner's Lemma for $1$-simplices.

\begin{Lem}\label{MScaleAndLZero}
Let $M$ be a positive number and $X$ be an
$M$-scale connected metric space. If $L^0({\cal {U}},X) > M$ for some
 cover ${\cal {U}}$ of $X$, then ${\cal {U}}$ contains
$X$ as an element.
\end{Lem}
\begin{pf} Suppose ${\cal {V}}$ is a refinement of ${\cal {U}}$
of multiplicity at most $1$ and Lebesque number bigger than $M$.
If $X$ is not an element of ${\cal {V}}$, then there are disjoint
non-empty elements $V_1,V_2\in {\cal {V}}$.
Pick a chain of points $x=x_1, x_2,\ldots, x_k=y$
such that $d_X(x_i,x_{i+1}) < M$ for all $i < k$ and $x\in V_1$, $y\in V_2$.
There is an index $j < k$ such that $x_j\in V_1$ and $x_{j+1}\notin V_1$.
The ball $B(x_{j+1},M)$ is contained in an element $W$ of ${\cal {V}}$
and intersects $V_1$. Therefore $W=V_1$, a contradiction.
\end{pf}

\section{The coarse category}

Let us introduce the coarse category in a way that explains
why two coarse functions are considered equivalent if their distance
is bounded.

\begin{Def} \label{CoarseFunDef}
Given a metric space $(X,d_X)$ and its two subsets $X_1$ and $X_2$
the notation $X_1\leq X_2$ means there is a positive number $R$
such that $X_1$ is contained in the ball $B(X_2,R)=\{x\in X \mid  dist(x,X_2) < R\}$.
\end{Def}

\begin{Prop}\label{CoarseFunIffRelation}
A function $f:X\to Y$ of metric spaces is coarse
if and only if it preserves the relation $\leq$ of sets. Thus, $X_1\leq X_2$ implies
$f(X_1)\leq f(X_2)$.
\end{Prop}
\begin{pf} Suppose $f:X\to Y$  preserves the relation $\leq$ of sets but not in the sense of Roe.
Therefore, for some $M > 0$ there is a sequence of points $x_n,y_n$
so that $d_X(x_n,y_n) < M$ for each $n$ but $d_Y(f(x_n),f(y_n))\to\infty$ as $n\to \infty$.
If $f(A)$ is bounded for some subsequence $A$ of $x_n$, then $f(B)$ is bounded for
the corresponding subsequence of $y_n$ (in view of $f(B)\leq f(A)$) contradicting
$d_Y(f(x_n),f(y_n))\to\infty$ as $n\to \infty$. Thus $f(x_n)\to\infty$ and $f(y_n)\to\infty$
as $n\to\infty$. By induction define a subsequence
$a_n$ of $\{x_n\}_{n\ge 1}$ and the corresponding subsequence
$b_n$ of $\{y_n\}_{n\ge 1}$ with the property that
$d_Y(f(a_k),f(b_i)) > k$ and $d_Y(f(b_k),f(a_i)) > k$ for all $k \ge i$.
Since $A=\{a_n\}_{n\ge 1}\leq B=\{b_n\}_{n\ge 1}$ one has $f(A)\leq f(B)$, a contradiction.

Suppose $f:X\to Y$ is coarse in the sense of Roe and $X_1\leq X_2$ in $X$.
Pick $R > 0$ so that $X_1\subset B(X_2,R)$ and choose $M > 0$ satisfying
$d_Y(f(x),f(y)) < M$ if $d_X(x,y) < R$ for all $x,y\in X$.
Given $x\in X_1$ pick $y\in X_2$ so that $d_X(x,y) < R$ since $d_Y(f(x),f(y)) < M$
one gets $f(X_1)\subset B(f(X_2),M)$. Thus $f(X_1)\leq f(X_2)$.
\end{pf}

Notice that $X_1\leq X_2$ for every bounded subset $X_1$ of $X$ provided $X_2\ne\emptyset$.
Also, $X_1\leq X_2$ implies $X_1$ is bounded provided $X_2$ is bounded.
Therefore $f(A)$ is bounded for every bounded subset $A$ of $X$ and every coarse function
$f:X\to Y$.

Given a function $f:X\to Y$ of metric spaces one can identify it with its graph
$\Gamma(f)\subset X\times Y$. Therefore it makes sense to ponder the meaning
of $\Gamma(f)\leq \Gamma(g)$ for $f,g:X\to Y$.

\begin{Prop}\label{FunEquivalentIffRelation}
Suppose $f,g:X\to Y$ are functions of metric spaces.
\begin{itemize}
\item[1.] If $g$ is coarse, then $\Gamma(f)\leq \Gamma(g)$ implies
that the distance $dist(f,g)$ between $f$ and $g$ is finite. In particular, $f$ 
is coarse.

\item[2.] If $dist(f,g)$ is finite,
then $\Gamma(f)\leq \Gamma(g)$.
\end{itemize}
\end{Prop}
\begin{pf} 1. Suppose the distance $dist(f,g)$ is not finite,
so there are points $x_n\in X$ with $d_Y(f(x_n),g(x_n) > n$ for all $n\ge 1$.
Let $R > 0$ be a number such that $B(\Gamma(g),R)$ contains $\Gamma(f)$.
For each $n$ pick $y_n\in X$ satisfying $d_X(x_n,y_n)+d_Y(f(x_n),g(y_n)) < R$.
There is $M > 0$ so that $d_Y(g(x_n),g(y_n)) < M$ for all $n\ge 1$ as $g$ is coarse.
Now, $d_Y(f(x_n),g(x_n)\leq d_Y(f(x_n),g(y_n))+d_Y(g(y_n),g(x_n)) < R+M$ for all $n\ge 1$,
a contradiction.

2. Notice $\Gamma(f)\subset B(\Gamma(g),dist(f,g))$.
\end{pf}

\begin{Def} \label{ForwardDistTransferDef}
Given a function $f:X\to Y$ of metric spaces we define
the {\it forward distance transfer} function $d_f:R_+\to R_+\cup\infty$
as the infimum of all functions $\alpha:R_+\to R_+\cup\infty$
with the property that $d_X(x,y)\leq t$ implies
$\alpha (t)\ge d_Y(f(x),f(y))$
for all $x,y\in X$.

The {\it reverse distance transfer} function $d^f:R_+\to R_+\cup\infty$
as the infimum of all functions $\alpha:R_+\to R_+\cup\infty$
with the property that $d_Y(f(x),f(y)))\leq t$
implies $d_X(x,y)\leq \alpha (t)$
for all $x,y\in X$.
\end{Def}

Notice that $f$ is coarse if and only if $d_f$ maps $R_+$ to $R_+$,
i.e. the values of $d_f$ are finite.
Also, $f$ is asymptotically Lipschitz if and only if
$d_f$ is bounded by a linear
function.

\begin{Prop}\label{FunEquivViaPreimages}
If $f,g:X\to Y$ are two coarsely proper coarse functions,
then the following conditions are equivalent:
\begin{itemize}
\item[1.] $dist(f,g)$ is finite.

\item[2.] For every coarse family ${\cal {U}}=\{U_s\}_{s\in S}$ in $Y$
the family $\{f^{-1}(U_s)\cap g^{-1}(U_s)\}_{s\in S}$ is coarse.
\end{itemize}
\end{Prop}
\begin{pf} 1$\implies$2. Let $dist(f,g) < M$. Consider
${\cal {V}}=\{B(U_s,-M)\}_{s\in S}$. It is a coarse family,
so $f^{-1}({\cal {V}})$ is coarse by \ref{ProperCoarseAndPtInv}. Notice
$f^{-1}(B(U_s,-M))\subset f^{-1}(U_s)\cap g^{-1}(U_s)$ for all $s\in S$
which is sufficient to establish coarseness of $\{f^{-1}(U_s)\cap g^{-1}(U_s)\}_{s\in S}$.

2$\implies$1. If $dist(f,g)$ is not finite, there is a sequence $x_n\to\infty$
such that $d_Y(f(x_n),g(x_n)) > n$ for all $n$. Put $A=\{x_n\}_{n\ge 1}$.
By \ref{CharOfFiniteCoarseFam}, the family ${\cal {U}}=\{Y\setminus f(A),Y\setminus g(A)\}$ is coarse.
However, $\{f^{-1}(U_s)\cap g^{-1}(U_s)\}_{s\in S}$ is not coarse
as it refines $\{X\setminus A\}$ which is not coarse.
\end{pf}

Our category is that of metric spaces and equivalence classes of coarse
functions. $f\sim g$ if $d_Y(f(x),g(x))$ is a bounded function of $x$.

Generalizing the concept of $A\leq B$ for subsets of a given metric space $X$,
we say $Y$ {\it coarsely dominates} $X$ (notation: $X\leq^{coa}_{rse} Y$)
if there are coarse functions $f:X\to Y$ and $g:Y\to X$ such that
$g\circ f$ is at a finite distance from $id_X$. 

\begin{Prop}\label{CoarseDomination}
Suppose $f:X\to Y$ and $g:Y\to X$ are coarse functions.
If $g\circ f$ is at a finite distance from $id_X$, then
 both $f:X\to f(X)$ and $g:f(X)\to X$ are coarsely proper
and $f\circ g$ is at finite distance from $id_{f(X)}$.
\end{Prop}
\begin{pf} Suppose $x_n\to \infty$. None of the
subsequences of $\{f(x_n)\}$ can be bounded as $g$ would
send it to a bounded subset of $X$. Thus $f(x_n)\to\infty$.
If $f(x_n)\to\infty$, then none of subsequences of $\{x_n\}$
is bounded. Therefore none of the subsequences of $\{g(f(x_n))\}$
is bounded and $g:f(X)\to X$ is coarsely proper.
If $d_X(g(f(x)),x) < M$ for all $x\in X$, then $d_Y(f(g(f(x))),f(x))\leq d_f(M)$
and $f\circ g$ is at finite distance from $id_{f(X)}$.
\end{pf}

\begin{Prop}\label{CharOfCoarseIso}
A surjective coarse function $f:X\to Y$ of metric spaces 
is a coarse isomorphism if and only if the reverse distance transfer
function $d^f$ is finite.
\end{Prop}
\begin{pf} If there is a coarse function $g:Y\to X$ such that
$g\circ f$ is at finite distance $M$ to $id_X$, then $d^f(a)\leq d_g(a)+2M$
is finite.

Assume $d^f$ is finite and pick a right inverse $g:Y\to X$.
Notice $d_X(g(x),g(y))\leq d^f(d_Y(x,y))$, so $g$ is coarse.
\end{pf}

\section{Coarse dimensions}

\begin{Def} \label{CoarseDimDef}
The {\it coarse dimension} $\dim^{coa}_{rse}(X)$ 
(respectively, {\it the major coarse dimension} $\dim^{COA}_{RSE}(X)$)
is the smallest integer $n$
such that any finite coarse family in $X$ (respectively, any coarse family in $X$)
has a coarse refinement with multiplicity at most $n+1$.
\end{Def}

\begin{Rem} \label{RemAboutCoarseDim}
Using Proposition 4.4 on p.1104 in \cite{Dr$_1$}
(notice that the words \lq uniformly bounded\rq \ are erroneously
inserted there) one can show that, for proper metric spaces $X$,
the major coarse dimension of $X$ coincides with the asymptotic dimension
of Dranishnikov. In view of \ref{CoarseDimAndHigsonCorona}, our coarse
dimension and Dranishnikov coarse dimension are identical.
\end{Rem}

Given a coarse family ${\cal {U}}=\{U_s\}_{s\in S}$ in a subset $A$ of $X$
one can extend it to a coarse family ${\cal {U}}'=\{U_s\cup(X\setminus A)\}_{s\in S}$
in $X$. Notice that ${\cal {V}}\cap A$ is a coarse refinement of ${\cal {U}}$
for any coarse refinement ${\cal {V}}$ of ${\cal {U}}'$.
Therefore the following holds.
\begin{Cor}\label{CoarseDimOfSubsets}
If $A$ is a subset of a metric space $X$, then
$\dim^{coa}_{rse}(A)\leq \dim^{coa}_{rse}(X)$
and $\dim^{COA}_{RSE}(A)\leq \dim^{COA}_{RSE}(X)$.
\end{Cor}

\begin{Prop}\label{CoarseDimAndDomination}
If $Y$ coarsely dominates $X$, then
$\dim^{coa}_{rse}(X)\leq \dim^{coa}_{rse}(Y)$
and $\dim^{COA}_{RSE}(X)\leq \dim^{COA}_{RSE}(Y)$.
\end{Prop}
\begin{pf} The proof is almost the same for both dimensions.
Suppose ${\cal {U}}$ is a coarse family in $X$ and $f:X\to Y$,
$g:Y\to X$ are coarse functions such that there is $M > 0$ with
$d_X(x,g(f(x))) < M$ for all $x\in X$. Replacing $Y$ by $f(X)$
 we may assume $f$ is onto
and both $f$ and $g$ are coarsely proper (see \ref{CoarseDomination}).
The idea of the proof is to refine $g^{-1}({\cal {U}})$ by ${\cal {V}}$
and then refine $f^{-1}({\cal {V}})$ to obtain a desired refinement
${\cal {W}}$ of ${\cal {U}}$ of multiplicity at most $n+1$, where $n$ is the 
dimension of $Y$.
Consider ${\cal {U}}'=\{B(U_s,-M)\}_{s\in S}$. It is a coarse family in $X$,
so $\{g^{-1}(B(U_s,-M))\}_{s\in S}$ is coarse and it has
a coarse shrinking ${\cal {V}}=\{V_s\}_{s\in S}$ of multiplicity at most $n+1$.
Suppose $x\in f^{-1}(V_s)\setminus U_s$. Since $d_X(x,g(f(x))) < M$,
$g(f(x))\notin B(U_s,-M)$. However, $f(x)\in V_s\subset g^{-1}(B(U_s,-M))$,
a contradiction.
\end{pf}

\begin{Def} \label{AsympDimDef}
The {\it minor asymptotic dimension} $ad(X)$ 
(respectively, {\it the asymptotic dimension} $asdim(X)$)
is the smallest integer $n$
such that the function ${\cal {U}}\to L^n({\cal {U}},X)$
is coarsely proper on the space of finite covers (respectively, arbitrary covers)
${\cal {U}}$ of $X$.
\end{Def}

Let us show that our definition of asymptotic dimension
is equivalent to that of Gromov.

\begin{Prop}\label{AsympDimEquivToGromov}
$asdim(X)\leq n$ if and only if for each $M > 0$ there is
a uniformly bounded family ${\cal {U}}$ in $X$ of Lebesque number at 
least $M$ and multiplicity at most $n+1$.
\end{Prop}
\begin{pf} If $asdim(X)\leq n$ as in \ref{AsympDimDef} and $M > 0$,
then there is $N > 0$ such that every cover ${\cal {V}}$ of $X$
satisfying $L({\cal {V}},X)\ge N$ has a refinement ${\cal {U}}$
of multiplicity at most $n+1$ and Lebesque number at least $M$.
Pick ${\cal {V}}$ to be the cover of $X$ by balls of radius $N$.
The resulting ${\cal {U}}$ is uniformly bounded.

\par Suppose for each $M > 0$ there is a uniformly bounded
family ${\cal {U}}^M$ of multiplicity at most $n+1$
and Lebesque number at least $M$. Let $\alpha(M)$ be the supremum of diameters
of elements of ${\cal {U}}^M$.
Given any family ${\cal {V}}$ of Lebesque number at least $\alpha(M)+1$,
${\cal {U}}^M$ is a refinement of of ${\cal {V}}$ which proves
that the function ${\cal {V}}\to L^n({\cal {V}},X)$
is coarsely proper on the space of all covers
${\cal {V}}$ of $X$.
\end{pf}

Quite often it is useful to have even stronger conditions imposed
on covers appearing in \ref{AsympDimEquivToGromov}.
\begin{Prop}[Gromov]\label{GromovOstrand} If 
Gromov asymptotic dimension $asdim(X)$ does not exceed $n$,
then for any $M, N > 0$ there exist uniformly bounded families
${\cal {U}}^i$, $1\leq i\leq n+1$, such that each ${\cal {U}}^i$
is $N$-disjoint and ${\cal {U}}=\bigcup\limits_{i=1}^{n+1}{\cal {U}}^i$
is of Lebesque number at least $M$.
\end{Prop}
\begin{pf}
Consider a uniformly bounded family ${\cal {V}}=\{V_s\}_{s\in s}$ of multiplicity at most $n+1$
and Lebesque number at least $2(n+1)\cdot (M+N)$.
Lemma \ref{OstrandPropertyForCovers} says it can be refined by ${\cal {V}}'=\bigcup\limits_{i=1}^{n+1}{\cal {V}}^i$
such that $L_{{\cal {V}}}(x)\ge L_{{\cal {U}}}(x)/(2n+2)\ge M+N$ for each $x\in X$
and each ${\cal {V}}^i$ consists of disjoint sets.
Define ${\cal {U}}_i$ as $\{B(W,-N)\}$, $W\in {\cal {V}}^i$.
\end{pf}

Let us characterize spaces of asymptotic dimension $0$.
\begin{Prop}\label{CharOfZeroAsympDim}
$asdim (X)>0$ if and only if there exist a number
$M>0$ and a coarsely proper sequence $\{(x_n,y_n)\}_{n=1}^\infty$ of pairs
of points in $X$ such that $dist(x_n,y_n)\to\infty$ and the points
$x_n$ and $y_n$ can be $M$-scale connected in $X\setminus
B(x_0,n)$.
\end{Prop}
\begin{pf} If $asdim(X)=0$, then for any $M>0$ there exists an
$M$-disjoint cover of $X$ by uniformly bounded sets. Therefore,
the distance between two points $x$ and $y$ which can be $M$-scale
connected in $X$ is uniformly bounded.
\par
Suppose $asdim (X)>0$. Let $n$ be a positive integer and $x_0$ be
the base point in $X$. There is $L>0$ such that $X$ does not have
a uniformly bounded cover of Lebesque number bigger than $L$ and
multiplicity 1. Define an equivalence relation on $X\setminus
B(x_0,n)$ by saying $x\sim y$ iff $x$ and $y$ can be $2L$-scale
connected in $X\setminus B(x_0,n)$. The cover of $X$ by the
equivalence classes has Lebesque number at least $2L$, therefore
these classes are not uniformly bounded by the choice of $L$.
Thus, there exist points $x_n$ and $y_n$ which can be $2L$-scale
connected in $X\setminus B(x_0,n)$ such that $dist(x_n,y_n)$ is
arbitrarily large.
\end{pf}

\begin{Prop}\label{AsympDimAndDomination}
If $Y$ coarsely dominates $X$, then
$asdim(X)\leq asdim(Y)$
and $ad(X)\leq ad(Y)$.
\end{Prop}
\begin{pf} The proof is almost the same for both dimensions.
Suppose ${\cal {U}}$ is a coarse family in $X$ and $f:X\to Y$,
$g:Y\to X$ are coarse functions such that there is $M > 0$ with
$d_X(x,g(f(x))) < M$ for all $x\in X$. By replacing
$Y$ with $f(X)$ we may assume $f$ is onto
and both $f$ and $g$ are coarsely proper (see \ref{CoarseDomination}).
The idea of the proof is to refine $g^{-1}({\cal {U}})$ by ${\cal {V}}$
and then refine $f^{-1}({\cal {V}})$ to obtain a desired refinement
${\cal {W}}$ of ${\cal {U}}$ of multiplicity at most $n+1$, where $n$ is the 
dimension of $X$.
Take a coarsely proper function $\alpha :R_+\to R_+$
with the property that any finite cover (respectively, arbitrary cover)
${\cal {U}}$ of $Y$ satisfying $L({\cal {U}},Y)\ge \alpha (t)$
has a refinement ${\cal {V}}$ of multiplicity at most $n+1$
so that $L({\cal {V}},Y)\ge t$.

Given $t > 0$ pick $\beta (t)$ so that $L^g(\beta (t)) > \alpha (t)$ 
(see \ref{CoarseFunctionDef}).
Assume $L({\cal {U}}) > M+\beta (t)$.
Consider ${\cal {U}}'=\{B(U_s,-M)\}_{s\in S}$. 
$L({\cal {U}}') > \beta  (t)$,
so $\{g^{-1}(B(U_s,-M))\}_{s\in S}$ is of Lebesque number at least $\alpha (t)$ and it has
a shrinking ${\cal {V}}=\{V_s\}_{s\in S}$ of multiplicity at most $n+1$
and $L({\cal {V}})\ge t$.
Suppose $x\in f^{-1}(V_s)\setminus U_s$. Since $d_X(x,g(f(x))) < M$,
$g(f(x))\notin B(U_s,-M)$. However, $f(x)\in V_s\subset g^{-1}(B(U_s,-M))$,
a contradiction.
\end{pf}

\begin{Thm}\label{MajorCoarseAndAsympDim} The major coarse dimension of $X$
does not exceed the asymptotic dimension of $X$.
\end{Thm}
\begin{pf}
Suppose $asdim(X)=n < \infty$ and ${\cal {U}}=\{U_s\}_{s\in S}$
is a coarse family in $X$. By Lemma \ref{ProperRefinementsOfCoarse} we may assume $U$ is coarsely proper.
By induction on $k$ find a sequence of numbers $M_0=1, M_1, M_2,\ldots$,
and covers ${\cal {V}}^k=\{V_t\}_{t\in T(k)}$, $k\ge 1$, 
of multiplicity at most $n+1$ and satisfying the following conditions:

a. $L({\cal {V}}^k,X)\ge M_{k-1}$ for $k\ge 1$.

b. The diameter of each element of ${\cal {V}}^k$ is smaller than $M_k$.

c. The family $\{B(x,M_{k-1}) \mid  d(x,x_0)\ge M_k\}$ refines ${\cal {U}}$ for each $k\ge 1$.

d. $M_{k+1} > 2M_k$ for all $k\ge 1$.

Find functions $j(k):T(k)\to T(k+1)$ so that $V_t\subset V_{j(k)(t)}$.
Denote $\{x: M_{k}\leq d(x,x_0) < M_{k+1}\}$ by $A_k$.
 Given $t\in T(k)$ so that $V_t$ is contained
in some element of ${\cal {U}}$ define $\alpha (t)\in S$ by looking
at the sequence $V_t\subset V_{j(k)(t)}\subset\ldots$, picking the latest
element contained in some $U_s$ and setting $\alpha (t)=s$ (it is possible
each element of the sequence is contained in some $U_s$ in which case
all of them are contained in some $U_s$ and that $s$ is picked as $\alpha(t)$).
Define $W_s$ as follows:
it is the union of non-empty sets of the form
$V_t\cap A_k$
so that $V_t\in{\cal {V}}^{k-1}$ and $\alpha (t)=s$.
Notice that $m({\cal {W}})\leq n+1$ as in the annulus $A_k$
the family
${\cal {W}}$ is obtained from ${\cal {V}}^{k-1}$ by assembling some of its elements together. 

We plan to show ${\cal {W}}$ is coarse by proving that if
$M_{k}\leq d(x,x_0) < M_{k+1}$, then $B(x,M_{k-3})$ is contained
in some $W_s$. Indeed, there is $t\in T(k-2)$ so that $B(x,M_{k-3})\subset V_t$.
Put $r=j(k-2)(t)$ and $u=j(k-1)(r)$. 
Points of $B(x,M_{k-3})$ can belong to only two of the following three annuli:
$A_{k-1}$, $A_k$,
and $A_{k+1}$.
If $z\in B(x,M_{k-3})\cap A_{k+1}$,
then $V_u\subset B(z,M_k)\subset U_s$ for some $s\in S$.
We might as well put $s=\alpha(t)=\alpha(u)=\alpha(r)$.
In this case $B(x,M_{k-3})\subset W_s$.
If $B(x,M_{k-3})$ misses the last annulus, then only $\alpha(r)$
is definitely defined ($\alpha(u)$ may not exist) and $\alpha(t)=\alpha(r)$.
Now, $B(x,M_{k-3})\subset W_s$, where $s=\alpha(r)$.
\end{pf}

\begin{Rem}\label{RemOnMajorAndAsdim} 
\ref{MajorCoarseAndAsympDim} generalizes Proposition 4.5 on p.1105 of 
\cite{Dr$_1$}.
\end{Rem}

\section{The large scale dimension}

In this section we prove that any dimension of $X$ (asymptotic, major coarse, or
minor asymptotic), if finite, equals
 the coarse dimension of $X$. That corresponds to results of Dranishnikov \cite{Dr$_1$}
that $asdim(X)$ or $asdim^\ast(X)$, if finite,
are equal to the dimension of the Higson corona of any proper metric space $X$.
Our proofs are direct and become simpler by introducing a new dimension,
the {\it large scale dimension} of $X$. That dimension turns out to be identical
with the coarse dimension.

\begin{Def}\label{LargeScaleDimDef} The {\it large scale dimension} $\dim^{large}_{scale}(X)$ of $X$
is the smallest integer $n$ such that $A\to L^{n}({\cal {U}},A)$
is a coarsely proper function on the set of bounded subsets of $X$
for all finite coarse families ${\cal {U}}$ in $X$.
\end{Def}

Notice $\dim^{large}_{scale}(X)=-1$ if $X$ is bounded.

Obviously, $\dim^{large}_{scale}(X)\ge \dim^{large}_{scale}(A)$
for any subset $A$ of $X$.

\begin{Prop}\label{LargeScaleDimIsTheSmallest}
 $ad(X)\ge \dim^{large}_{scale}(X)$ and $\dim^{coa}_{rse}(X)\ge \dim^{large}_{scale}(X)$.
\end{Prop}
\begin{pf} The inequality $\dim^{coa}_{rse}(X)\ge \dim^{large}_{scale}(X)$ is almost obvious.
Indeed, given $n=\dim^{coa}_{rse}(X)$ and given a coarse family ${\cal {U}}$ in $X$
consisting of $m$ elements one has a coarse refinement ${\cal {V}}$ of ${\cal {U}}$
such that the multiplicity $m({\cal {V}})$ is at most $n+1$.
In that case 
$$L^{n}({\cal {U}},A)\ge L({\cal {V}},A)\ge \inf\limits_{a\in A}L_{{\cal {V}}}(a)$$
and is a coarsely proper function of $A$.

Suppose  $ad(X)=n$ and ${\cal {U}}$ is a coarse cover of $X$ consisting
of $m$ elements. Given $t > 0$ find a bounded subset $U$ of $X$
such that ${\cal {U}}|_ {(X\setminus U)}$ has a refinement ${\cal {V}}$
of multiplicity at most $n+1$ and Lebesque number at least $t$.
For any bounded subset $A$ of $X\setminus U$, $L^{n}({\cal {U}},A)\ge L({\cal {V}},A)\ge t$
which proves $\dim^{large}_{scale}(X)\leq n$.
\end{pf}

As shown in \cite{DKU}, the asymptotic dimension of $R^n$
is at most $n$ (see p.793). For the convenience of the reader
let us reword the argument from \cite{DKU} as follows:
Given $M > 0$ consider the triangulation on the unit $n$-cube $I^n$ obtained
by starring at the center of each face. It is invariant under symmetries
of $I^n$ and the cover of $I^n$ by stars of vertices has a positive Lebesque number
$k$ and is of multiplicity at most $n+1$. Rescale $I^n$ by the factor of $M/k$
and extend its triangulation over the whole $R^n$ by reflections.
The cover of $R^n$ by stars of vertices has Lebesque number at least $M$ and
is of multiplicity at most $n+1$.

Let us show how to use the large scale dimension to estimate asymptotic
dimension from below.

\begin{Prop}\label{DimOfRn} $\dim^{large}_{scale}(R^n)\ge n$.
\end{Prop}
\begin{pf} Since $\dim(I^n)=n$, there is a finite open cover ${\cal {U}}$ of $I^n$
with no open refinement of multiplicity at most $n$.
Let $I^n_k\subset R^n$ be a copy of $I^n$ enlarged $k$ times with the corresponding
cover ${\cal {U}}^k$. We request $I^n_k\to\infty$ so that ${\cal {V}}$ obtained
by adding the corresponding elements of ${\cal {U}}^k$ is a finite coarse family
on $A=\bigcup\limits_{k=1}^\infty I^n_k$.
Notice $L^{n-1}({\cal {V}},I^n_k)=0$ for all $k$. Thus
$\dim^{large}_{scale}(R^n)\ge n$.
\end{pf}

\begin{Prop}\label{ZeroDimCase} If $\dim^{large}_{scale}(X)=0$, then $asdim(X)=0$
and $\dim^{COA}_{RSE}(X)=0$.
\end{Prop}
\begin{pf}
It suffices to show $asdim(X)=0$ (see  \ref{MajorCoarseAndAsympDim}).
Suppose $asdim(X) > 0$. By \ref{CharOfZeroAsympDim} there exist a number $M>0$
and a coarsely proper sequence $\{(x_n,y_n)\}_{n=1}^\infty$ of pairs of
points in $X$ such that $dist(x_n,y_n)\to\infty$ and the points
$x_n$ and $y_n$ can be $M$-scale connected in $X\setminus
B(x_0,n)$ by a chain $P_n$. Consider a coarse family ${\cal {U}}$
consisting of two sets: $X\setminus \bigcup\limits_{n=1}^\infty\{x_n\}$ and
$X\setminus \bigcup\limits_{n=1}^\infty\{y_n\}$.

Since $C\to L^0({\cal {U}},C)$ is a coarsely proper function, there is a chain
$P_n$ such that $L^0({\cal {U}},P_n)>M$. This contradicts \ref{MScaleAndLZero} since
$P_n$ is $M$-scale connected and the cover ${\cal {U}}$ is non-trivial on
$P_n$.
\end{pf}

\begin{Def}\label{StandardPU} Given a point-finite family ${\cal {U}}=\{U_s\}_{s\in S}$
in $X$
(that means each point of $X$ belongs to at most finitely many elements of ${\cal {U}}$)
by {\it the canonical partition of unity} of ${\cal {U}}$ we mean the family
of functions $\{f_s/f\}_{s\in S}$, where $f_s(x)=dist(x,X\setminus U_s)$
and $f(x)=\sum\limits_{s\in S} f_s(x)$.
If $T$ is a subset of $S$, then $X_T$ is defined to be
$\{x\in X \mid  \sum\limits_{s\in T} f_s(x)/f(x)=1\}$ and
by $\partial X_T$ we mean the set of all $x\in X_T$ such that $f_s(x)=0$ for some $s\in T$.
\end{Def}

Notice that $f(x) > 0$ for all $x\in X$ such that $L_{{\cal {U}}}(x) > 0$
and $f$ is a Lipschitz function if ${\cal {U}}$ is of finite multiplicity.

\begin{Lem}\label{LargScaleDimAndRefinements} If the large scale dimension of $X$ is at most $n$,
then any coarse family ${\cal {U}}$ in $X$ of finite multiplicity $m$
has a coarse refinement ${\cal {V}}$ of multiplicity at most $n+1$.
\end{Lem}
\begin{pf} Suppose ${\cal {U}}$ exists with no coarse
refinement of multiplicity at most $n+1$. Using \ref{ProperRefinementsOfCoarse} we reduce
the general case to that of ${\cal {U}}=\{U_s\}_{s\in S}$ consisting
of bounded sets so that for any sequence $x_k\to\infty$
the conditions $x_k\in U_{s(k)}\in{\cal {U}}$ imply $U_{s(k)}\to\infty$.
For induction on $m-n$ it suffices to assume the multiplicity of
${\cal {U}}$ is $n+2$.
\par

 Pick a coarse shrinking ${\cal {W}}=\{W_s\}_{s\in S}$ (see \ref{CoarseShrinkingALaParacompact})
so that given $M > 0$ there is a bounded subset $A$ of $X$
with the property that, for $x\in X\setminus A$,
$B(x,M)\cap W_s\ne\emptyset$ implies $B(x,M)\subset U_s$.
Consider the canonical partition of unity $f$ of ${\cal {W}}$.
Given a set $T$ in $S$ consisting of $n+2$ elements pick a shrinking ${\cal {W}}^T$ of
${\cal {W}}|_ {X_T}$ of
order at most $n+1$ and the Lebesque number at least half the maximum 
$L^n({\cal {W}}^T,X_T)$ possible
(if the maximum is infinity we pick a shrinking of Lebesque number twice
the size of $X_T$). We can add $W_s\cap \partial X_T$ to $W^T_s$ without
increasing the order of $W^T$ beyond $n+1$ (obviously, the Lebesque number does not decrease).
By pasting those shrinkings for all $T$
 one gets a refinement ${\cal {V}}$ of ${\cal {W}}$
on $X\setminus A$ for some bounded subset $A$ of $X$ of multiplicity at most $n+1$.
Therefore ${\cal {V}}$ cannot be coarse and there is $M > 0$ and a sequence
of points $x_k\to\infty$ such that none of $B(x_k,M)$ is contained
in an element of ${\cal {V}}$. In particular $B(x_k,M)$ is not contained in the
$n$-skeleton of $X$ (the points where the order of $f$
is at most $n+1$) for large $k$.

Pick sets $T(k)$ so that $X_{T(k)}\setminus\partial X_{T(k)}$
contains an element $y_k\in B(x_k,M)$. For large $k$, $B(x_k,M)$ intersecting
$W_s$ implies $B(x_k,M)\subset U_s$. Therefore the set $T$ of $s\in S$
so that $B(x_k,M)$ intersects
$W_s$ is of cardinality at most $n+2$ and
 $B(x_k,M)\subset X_{T(k)}$.
For large $k$ the cover ${\cal {W}}|_ {X_{T(k)}}$ has a refinement of order
at most $n+1$ and Lebesque number at least $3M$. Therefore, $B(x_k,M)$
is contained in a single element of ${\cal {V}}$, a contradiction.
\end{pf}

\begin{Cor}\label{CoarseDimEqualsLgScale} The coarse dimension of $X$
 equals the large scale dimension of $X$.
\end{Cor}

\begin{Cor}\label{MajorCoarseEqualsCoarseIfFinite} If the major coarse dimension of $X$ 
is finite, then it equals the large scale dimension of $X$.
\end{Cor}

\begin{Thm} \label{AsympDimEqualsCoarseIfFinite}
If the asymptotic dimension (respectively, the minor asymptotic dimension)
 of unbounded $X$ is finite, then it equals
the large scale dimension of $X$.
\end{Thm}
\begin{pf}
Suppose $asdim(X)=n$ (respectively, $ad(X)=n$) and $\dim^{large}_{scale}(X)< n$.
Notice $n > 0$ as $\dim^{large}_{scale}(X)< 0$ is possible only for bounded $X$.
Therefore there is $M > 0$ and a sequence of covers
(respectively, finite covers) ${\cal {U}}^k$ indexed by sets $S(k)$ of Lebesque number at least $k+3M$
and multiplicity at most $n+1$
so that no refinement of ${\cal {U}}^k$
of multiplicity $n$ has Lebesque number bigger than $M$.
Augment each ${\cal {U}}^k$ by shrinking it
to the family $B(U,-M)$, $U\in {\cal {U}}^k$. Let $f^k$ be the canonical
partition of unity of that augmentation.
\par
Notice that for any $k$ and any $x\in X$
there is a subset $T$ of $S(k)$ consisting of at most $(n+1)$ elements
so that $B(x,M)\subset X_T$.
We are going to show that for every $k$ there
is $N > 0$ such that for any $R > N$ there is $T(k)\subset S(k)$
consisting of at most $(n+1)$ elements with
$X_{T(k)}\subset X\setminus B(x_0,R)$, $x_0$ a fixed point in $X$,
so that $Carr(f^k|_ {X_{T(k)}})$ does not admit a refinement
of multiplicity at most $n$ and Lebesque number bigger than $M$.
\par
Suppose that, for some $k$ and $R > 0$, all $Carr(f^k|_ {X_T})$
so that $X_{T}\subset X\setminus B(x_0,R)$ do admit a  
refinement ${\cal {V}}(T)$ of multiplicity at most $n$ and Lebesque number bigger than $M$.
By converting those refinements to shrinkings and pasting
 one gets a refinement ${\cal {V}}$ of ${\cal {U}}^k$
on $X\setminus U$ for some bounded subset $U$ of $X$ of multiplicity at most $n$
and Lebesque number bigger than $M$. More precisely, for each $T\subset T(k)$
so that $X_{T}\subset X\setminus B(x_0,R)$, we pick
a shrinking $\{V^T_t\}_{t\in T}$ of $Carr(f^k|_ {X_T})$ of multiplicity at most $n$
and Lebesque number bigger than $M$. If $T$ contains at most $n$ elements,
that shrinking is picked to be exactly $Carr(f^k|_ {X_T})$ as the multiplicity
is at most $n$ in such case. ${\cal {V}}$ is a shrinking of ${\cal {U}}^k|_{(X\setminus U)}$, $U$ being the union of $X_T$ that are not contained in
$X\setminus B(x_0,R)$,
and $V_s$, $s\in S(k)$, is defined as the union of all $V^T_s$ with $s\in T$. The reason
${\cal {V}}$ has Lebesque number at least $M$ is that for any $x\in X$
there is a subset $T$ of $S(k)$ consisting of at most $(n+1)$ elements
so that $B(x,M)\subset X_T$. 
\par
Now, the cover consisting of the union of
$B(U,2M)$ and all the elements of ${\cal {V}}$ intersecting $B(U,2M)$
and of all elements of ${\cal {V}}$ that do not intersect $B(U,2M)$ is uniformly bounded, of multiplicity at most $n$ (recall $n > 0$),
and of Lebesque number bigger than $M$, a contradiction.
\par
Construct by induction a sequence of sets $T(i)\subset S(i)$ 
with $X_{T(i)}$ being mutually disjoint and tending to infinity
so that $Carr(f^i|_ {X_{T(i)}})$ does not have a refinement of multiplicity
at most $n$ and Lebesque number bigger than $M$.
Paste all those carriers according to their index within each set $T(i)$
and get a coarse cover on a subset $A$ of $X$ that does not admit
a refinement of multiplicity at most $n$ and Lebesque number bigger than $M$
on infinitely many $X_{T(i)}$, a contradiction.
\end{pf}

\section{Slowly oscillating functions }

\begin{Def}\label{SlowlyOscFunDef}
A function $f:X\to Y$ is slowly oscillating if
$f^{-1}({\cal {U}})$ is coarse for every cover ${\cal {U}}$ of $Y$
of positive Lebesque number.
\end{Def}

\begin{Def} \label{OscDef}
Given a function $f:X\to Y$ of metric spaces one defines
its {\it oscillation function} $Osc(f,M):X\to R_+\cup\infty$
for every $M > 0$ by declaring $Osc(f,M)(a)$ to be the supremum
of $d_Y(f(x),f(a))$ over all $x\in B(a,M)$.
\end{Def}

\begin{Prop}\label{OurOscDefSameAsClassical}
$f$ is slowly oscillating if and only if $Osc(f,M)(x)\to 0$ as $x\to \infty$
for all $M > 0$.
\end{Prop}
\begin{pf} Suppose $Osc(f,M)(x)\to 0$ as $x\to \infty$
for all $M > 0$. Given a cover ${\cal {U}}$ of $Y$
of positive Lebesque number and given $x_n\to\infty$ in $X$
there is $N > 0$ such that each $f(B(x_n,M))$ is of diameter smaller
that $L({\cal {U}},Y)$ for $n > N$.
Therefore $B(x_n,M)$ is contained in an element of $f^{-1}({\cal {U}})$ 
and $f^{-1}({\cal {U}})$ is coarse.
\par Suppose $f^{-1}({\cal {U}})$ is coarse for every cover ${\cal {U}}$ of $Y$
of positive Lebesque number.
Given $x_n\to \infty$ in $X$ and given $M > 0$ such that
diameters of $f(B(x_n,M))$ are bigger than a fixed $\delta  > 0$,
consider ${\cal {U}}=\{B(y,\delta /2)\}_{y\in Y}$.
Since $f^{-1}({\cal {U}})$ is coarse, there is $N > 0$ such that for
all $n > N$ sets $B(x_n,M)$ are contained in an element of
$f^{-1}({\cal {U}})$. Therefore diameters of $f(B(x_n,M))$ are smaller than a $\delta$
for $n > N$, a contradiction.
\end{pf}

Our basic way of constructing slowly oscillating real-valued functions
is based on the following. 

\begin{Lem} \label{OscOfFraction}
Suppose $f,g:X\to R_+$ and $Osc(f,M),Osc(g,M) < \epsilon$ for some $\epsilon > 0$.
If $f(x)+g(x) > N$ for all $x\in X$, then $Osc(\frac{f}{f+g},M) < \frac{3\epsilon}{N}$.
\end{Lem}
\begin{pf}
Let $h=\frac{f}{f+g}$ and $a=\frac{3\epsilon}{N}$. If $h(x)-h(y) \ge a$
for some $x,y\in X$ satisfying $d_X(x,y)< M$, then
$\frac{f(x)}{f(x)+g(x)}-\frac{f(x)-\epsilon}{f(x)+g(x)+2\cdot \epsilon} \ge a$ as well.
Since $\frac{f(x)}{f(x)+g(x)}-\frac{f(x)-\epsilon}{f(x)+g(x)+2\cdot \epsilon} =
\frac{f(x)\cdot 2\epsilon+\epsilon\cdot (f(x)+g(x))}{(f(x)+g(x))\cdot (f(x)+g(x)+2\epsilon)}
\leq \frac{3\epsilon}{f(x)+g(x)+2\epsilon}< a$,
we arrive at a contradiction.
\end{pf}

\begin{Cor} \label{FractionIsSlowOsc}
If $f$ and $g$ are coarse
functions from $X$ to $R_+$ such that $f+g$ is coarsely proper and positive, then
$f/(f+g)$ is slowly oscillating.
\end{Cor}

Here is a simple connection between oscillation and the Lebesque number.

\begin{Lem}\label{OscAndLebesqueNum}
If $\phi=\{\phi_s:X\to R_+\}_{s\in S}$ is a family of functions
with finite supremum $\sup(\phi)$
such that $Osc(\phi_s,M) < \frac{1}{2}\sup(\phi)$ for each $s\in S$,
then $L(\phi)\ge M$.
\end{Lem}
\begin{pf}
Given $a\in X$ find $s\in S$ so that $\phi_s(a) > \frac{1}{2}\sup(\phi)(a)$.
If $d_X(x,a) < M$,
then $|\phi_s(x)-\phi_s(a)| < \frac{1}{2}\sup(\phi)(a)$, so $\phi_s(x)$
cannot be $0$ thus affirming $B(a,M)\subset \phi_s^{-1}(0,\infty)$.
\end{pf}

 A partition of unity $\phi=\{\phi_s:X\to R_+\}_{s\in S}$
 is called {\it slowly oscillating} if the corresponding
function $\phi:X\to l^1_S$ is slowly oscillating.

$\phi$ is called {\it equi-slowly oscillating} if 
 the oscillation of all $\phi_s$ is synchronized in the following way:
for every $M > 0$ and every $\epsilon > 0$ there is a bounded subset $U$
of $X$ such that $Osc(\phi_s,M)(x) < \epsilon$ for all $x\in X\setminus U$
and all $s\in S$.
Obviously, every finite partition of unity into slowly oscillating functions
is globally slowly oscillating and is equi-slowly oscillating.
Also, every slowly oscillating partition of unity is equi-slowly oscillating.

\begin{Lem}\label{SlowOscPUAndEquiosc}
If $\phi=\{\phi_s:X\to R_+\}_{s\in S}$ is a partition of unity
of finite multiplicity $m$, then $\phi$ is slowly oscillating if and only if
it is equi-slowly oscillating.
\end{Lem}
\begin{pf}
Given $M,\epsilon > 0$ we can find a bounded set $U$
such that $Osc(\phi_s,M) < \epsilon/(2m)$ for all $x\in X\setminus U$ and all $s\in S$.
If $a\in X\setminus U$ and $x\in B(a,M)$,
then the complement $F$ of set $T=\{s\in S \mid  \phi_s(x)+\phi_s(a)=0\}$ contains at most $2m$
elements. Since $|\phi(x)-\phi)(a)| = \sum\limits_{s\in F}|\phi_s(x)-\phi_s(a)|< |F|\cdot \epsilon/(2m)
\leq\epsilon$, $\phi$ is slowly oscillating.
\end{pf}

\begin{Lem}\label{EquiSlowOscAndCoarseCarriers}
If $\phi=\{\phi_s:X\to R_+\}_{s\in S}$ is an equi-slowly oscillating partition of unity
of finite multiplicity $m$, then its carrier family $Carr(\phi)$ is coarse.
\end{Lem}
\begin{pf}
Notice $\sup(\phi) \ge 1/m$.
Given $M > 0$ we can find a bounded set $U$
such that $Osc(\phi_s,M) < 1/(2m)$ for all $x\in X\setminus U$ and all $s\in S$.
By \ref{OscAndLebesqueNum}, $L(\phi|_{(X\setminus U)},X\setminus U) > M$ which proves $Carr(\phi)$ is coarse.
\end{pf}

\begin{Rem}\label{RemarkOnPUAndCoarse}
If one drops the assumption of $\phi$ being of finite multiplicity,
then the carrier family may not be coarse:
Take a cloud $C_n$ of $2^n+1$ points at location $2^n$ with mutual distances equal $1$.
For each $x\in X$ define $\phi_x$ as taking value $0$ at $x$ and all points not in its cloud.
For points $y\in Cloud(x)\setminus\{x\}$ we put $\phi_x(y)=2^{-n}$.
\end{Rem}

\begin{Cor}\label{FiniteCoarseAndSlowOscPUs}
If ${\cal {U}}=\{U_s\}_{s\in S}$ is a cover of $X$ of finite multiplicity,
then the following conditions are equivalent:
\begin{itemize}
\item[1.] ${\cal {U}}$ is coarse.

\item[2.] There is a continuous slowly oscillating partition of unity
$\phi=\{\phi_s\}_{s\in S}$ on $X\setminus A$
for some bounded subset $A$ of $X$ such that $Carr(\phi_s)\subset U_s$
for each $s\in S$.
\item[3.] There is a slowly oscillating partition of unity
$\phi=\{\phi_s\}_{s\in S}$ on $X\setminus A$
for some bounded subset $A$ of $X$ such that $Carr(\phi_s)\subset U_s$
for each $s\in S$.
\end{itemize}
\end{Cor}
\begin{pf} 1$\implies$2. Define $f(x)=\sum\limits_{s\in S} dist(x,X\setminus U_s)$ and
$f_s(x)=dist(x,X\setminus U_s)$.
Notice that $f$ is a coarsely proper Lipschitz function and
 \ref{FractionIsSlowOsc} says that $\{f_s/f\}_{s\in S}$
is an equi-slowly oscillating partition of unity on $X\setminus A$,
where $A$ is the zero-set of $f$. By  \ref{SlowOscPUAndEquiosc} it is
a slowly oscillating partition of unity.

2$\implies$3 is obvious.

3$\implies$1 follows from \ref{EquiSlowOscAndCoarseCarriers}.
\end{pf}

\section{Coarse dimension and Higson corona}

Given a metric space $X$ by the {\it Higson compactification} of $X$
we mean a compact Hausdorff space $h(X)$ containing $X$ as a dense subset with the property
that a bounded continuous function $f:X\to R_+$ extends over $h(X)$
if and only if $f$ is slowly oscillating.
If the metric on $X$ is proper and $X$ is locally compact,
then $X$ is open in $h(X)$ and the remainder $h(X)\setminus X$
is called the {\it Higson corona} of $X$ and denoted by $\nu(X)$.

A metric space $X$ is called {\it $\delta $-disjoint} for some $\delta  > 0$
if $d_X(x,y)\ge \delta $ for all $x\ne y$.

\begin{Thm} \label{CoarseDimOfHigsonComp}
If $X$ is a $\delta $-disjoint metric space for some $\delta > 0$,
then its coarse dimension equals the dimension of the Higson 
compactification of $X$.
\end{Thm}
\begin{pf} Suppose $\dim^{coa}_{rse}(X)=m < \infty$.
Given a finite open cover ${\cal {U}}=\{U_s\}_{s\in S}$ of the Higson compactification $h(X)$ of
$X$ we find a partition of unity $f=\{f_s\}_{s\in S}$ on $h(X)$
such that $cl(f_s^{-1}(0,1])\subset U_s$ for each $s\in S$ (see \cite{Dy}).
As $f|_ X$ is slowly oscillating (see \ref{SlowOscPUAndEquiosc}), the family
$\{f_s^{-1}(0,1]\cap X\}_{s\in S}$ is coarse in $X$ (see \ref{EquiSlowOscAndCoarseCarriers}).
By \ref{FiniteCoarseAndSlowOscPUs} there is a slowly oscillating partition of unity
$g=\{g_s\}_{s\in S}$ on $X$ whose multiplicity is at most $m+1$
and $g_s^{-1}(0,1]\subset f_s^{-1}(0,1]\cap X$ for each $s\in S$.
Extend each $g_s$ over $h(X)$ to $k_s:h(X)\to [0,1]$.
The resulting family $k=\{k_s\}_{s\in S}$ is a partition of unity on $h(X)$.
It remains to show $m(k)\leq m+1$ and $k_s^{-1}(0,1]\subset U_s$ for each $s\in S$.
If there is a point $x\in h(X)\setminus X$ such that $k_s(x) > 0$
for all $s\in T$, $T$ containing at least $m+2$ elements,
then the same would be true for some neighborhood $U_x$ of $x$ in $h(X)$.
Since $U_x\cap X\ne\emptyset$ one arrives at a contradiction with the fact
that $m(g)\leq m+1$. If $k_s^{-1}(0,1]$ is not a subset of $U_s$
for some $s\in S$, then there is $x\in h(X)\setminus X$ so that
$x\in k_s^{-1}(0,1]\setminus cl(f_s^{-1}(0,1])$.
That means there is a neighborhood $U_x$ of $x$ in $h(X)$ on which
$f_s$ is identically $0$. Hence $g_s|_ {(U_x\setminus X)}\equiv 0$
implying $k_s(x)=0$, a contradiction.
\end{pf}

\begin{Cor}\label{CoarseDimAndHigsonCorona}
If $X$ is a proper metric space, then the dimension
of its Higson corona equals the coarse dimension of $X$.
\end{Cor}
\begin{pf} Consider a maximal $1$-disjoint subset $A$ of $X$.
Notice $\dim^{coa}_{rse}(A)=\dim^{coa}_{rse}(X)$
and Higson coronas $\nu(A)$ and $\nu(X)$ for both $A$ and $X$ are identical.
Since $A$ is $1$-disjoint, $\dim^{coa}_{rse}(A)=\dim(h(A))=\dim(\nu(A))=\dim(\nu(X))$.
\end{pf}

\begin{Cor} \label{UnionThmForCoarse}
If $X=A\cup B$,
then the coarse dimension of $X$ equals maximum of the coarse dimensions
of $A$ and $B$.
\end{Cor}
\begin{pf} Let $m=\max(\dim^{coa}_{rse}(A),\dim^{coa}_{rse}(B))$.
By \ref{CoarseDimOfSubsets}, $\dim^{coa}_{rse}(X)\ge m$.
By switching to maximal $1$-disjoint subsets of $A$ and $B$, respectively,
we reduce the general case to that of $X$ being $1$-disjoint.
Consider the Higson compactification $h(X)$ of $X$. Notice $cl(A)$
is the Higson compactification of $A$ as any slowly oscillating and bounded function
$f:A\to R_+$ extends over $X$ to a bounded and slowly oscillating function.
The same is true for $B$. Since $h(X)=cl(A)\cup cl(B)$,
$\dim(h(X))=\max(\dim(cl(A)),\dim(cl(B)))=\max(\dim^{coa}_{rse}(A),\dim^{coa}_{rse}(B))=m$.
\end{pf}

We plan to extend \ref{UnionThmForCoarse} to other dimensions as well.
Our strategy is to show finiteness of the appropriate dimension
of $X$ first, then use \ref{UnionThmForCoarse} as well as the fact that
all other dimensions are equal to the coarse dimension of $X$
once they are finite (see \ref{MajorCoarseEqualsCoarseIfFinite} and
\ref{AsympDimEqualsCoarseIfFinite}).

\begin{Cor} \label{UnionThmForAsymptotic}
If $X=A\cup B$,
then the asymptotic dimension of $X$ equals maximum of the asymptotic dimensions
of $A$ and $B$.
\end{Cor}
\begin{pf} Let $m=\max(asdim(A),asdim(B))$.
Obviously $asdim(X)\ge m$.
Given $M > 0$ find uniformly bounded family
${\cal {U}}_A$ in $A$ covering $A$ and being the union of $m+1$ families, each of them $3M$-disjoint.
Similarly, find uniformly bounded and $3M$-disjoint family
${\cal {U}}_B$ in $B$ covering $B$ and
 being the union of $m+1$ families, each of them $3M$-disjoint.
Consider ${\cal {U}}={\cal {U}}_A\cup {\cal {U}}_B$ and let ${\cal {V}}=\{B(U,M)\}_{U\in {\cal {U}}}$.
Notice ${\cal {V}}$ is uniformly bounded in $X$, is of multiplicity
at most $2(m+1)$, and $L({\cal {V}},X)\ge M$.
Therefore $asdim(X)\leq 2m+1$ and (see
\ref{AsympDimEqualsCoarseIfFinite}) $asdim(X)=\dim^{coa}_{rse}(X)=m$.
\end{pf}
\begin{Rem}\label{RemarkAboutFiniteUnionThm}
\ref{UnionThmForAsymptotic}
was proved in \cite{BD$_1$} (see the Finite Union Theorem there) for $X$
being a proper metric space by using totally different methods.
\end{Rem}

\begin{Cor} \label{UnionThmForMajorCoarse}
If $X=A\cup B$,
then the major coarse dimension of $X$ equals maximum of the major coarse dimensions
of $A$ and $B$.
\end{Cor}
\begin{pf} Let $m=\max(\dim^{coa}_{rse}(A),\dim^{coa}_{rse}(B))$.
By \ref{CoarseDimOfSubsets}, $\dim^{COA}_{RSE}(X)\ge m$.
Given a coarse family ${\cal {U}}$ in $X$ put $f(x)=L_{{\cal {U}}}(x)$.
If $f(x)=\infty$ for some $X$, then ${\cal {U}}$ has a coarse refinement
of order at most $2$ (see \ref{CoversWithInfiniteLebesque}). Assume $f(x)<\infty$ for all $x\in X$.
Pick a coarse refinement $\{V_a\}_{a\in A}$ of multiplicity
at most $m+1$ of the family $\{B(a,f(a)/2)\}_{a\in A}$.
Pick a coarse refinement $\{V_b\}_{b\in B}$ of multiplicity
at most $m+1$ of the family $\{B(b,f(b)/2)\}_{b\in B}$.
If $V_a\ne\emptyset$ define $e(V_a)=\{x\in B(a,f(a)) \mid  dist(x,V_a) < dist(x,A\setminus V_a\}$.
Observe $\bigcap\limits_{a\in T}e(V_a)\ne\emptyset$ implies 
$\bigcap\limits_{a\in T}V_a\ne\emptyset$ for every finite subset $T$ of $S$.
Indeed, suppose $x\in \bigcap\limits_{a\in T}e(V_a)$ and find $\delta > 0$ such that
$dist(x,V_a)+\delta  < dist(x,A\setminus V_a\}$ for all $a\in T$.
Pick $y\in A$ so that $dist(x,A)+\delta > d(x,y)$.
If $y\in A\setminus V_a$ for some $a\in T$, then
$dist(x,A)+\delta \leq dist(x,V_a)+\delta < dist(x,A\setminus V_a)\leq d(x,y)$,
a contradiction.
Therefore the multiplicity of $\{e(V_a)\}_{a\in A}$ is at most $m+1$.
Do the same procedure for $B$ and produce $\{e(V_b)\}_{b\in B}$.
If $x\in a$, $M < f(a)$ and $B(x,M)\cap A\subset V_a$, then $B(x,M/2)\subset e(V_a)$.
Therefore $\{e(V_a)\}_{a\in A}\cup \{e(V_b)\}_{b\in B}$ is coarse in $X$,
refines ${\cal {U}}$, and is of multiplicity at most $2(m+1)$.
Thus $\dim^{COA}_{RSE}(X)\leq 2m+1$ and (see 
\ref{MajorCoarseEqualsCoarseIfFinite}) $\dim^{COA}_{RSE}(X)=\dim^{coa}_{rse}(X)=m$.
\end{pf}

\begin{Cor} \label{UnionThmForMinorAsymptotic}
If $X=A\cup B$,
then the minor asymptotic dimension of $X$ equals maximum of the minor asymptotic dimensions
of $A$ and $B$.
\end{Cor}
\begin{pf} Let $m=\max(ad(A),ad(B))$.
Obviously $ad(X)\ge m$.
Suppose $M > 0$ and find $N > 0$ such that $L^m({\cal {U}},A)> 2M$
for all finite covers ${\cal {U}}$ of $A$ satisfying $L({\cal {U}},A) > N$.
We can use the same $N$ and claim $L^m({\cal {U}},B)> 2M$
for all finite covers ${\cal {U}}$ of $B$ satisfying $L({\cal {U}},B) > N$.
Given a finite family ${\cal {U}}=\{U_s\}_{s\in S}$ in $X$ 
satisfying $L({\cal {U}},X) > M+N$, consider $\{B(U_s,-M)\}_{s\in S}$
and shrink it on $A$ to a family $\{V_s\}_{s\in S}$ of multiplicity
at most $m+1$ and Lebesque number at least $2M$.
Do the same for $B$ and shrink $\{B(U_s,-M)\}_{s\in S}$
 on $B$ to a family $\{W_s\}_{s\in S}$ of multiplicity
at most $m+1$ and Lebesque number at least $2M$.
If $V_s\ne\emptyset$ define $e(V_s)=\{x\in U_s \mid  dist(x,V_s) < dist(x,A\setminus V_s\}$.
Observe $\bigcap\limits_{s\in T}e(V_s)\ne\emptyset$ implies 
$\bigcap\limits_{s\in T}V_s\ne\emptyset$ for every finite subset $T$ of $S$
(see the proof of  \ref{UnionThmForMajorCoarse}).
Therefore the multiplicity of $\{e(V_s)\}_{s\in S}$ is at most $m+1$.
Do the same procedure for $B$ and produce $\{e(W_s)\}_{s\in S}$.
Obviously $\{e(V_s)\}_{s\in S}\cup \{e(V_s)\}_{s\in S}$ 
refines ${\cal {U}}$ and is of multiplicity at most $2(m+1)$.
If we show its Lebesque number is at least $M$ we will demonstrate
 $ad(X)\leq 2m+1$ and (see \ref{AsympDimEqualsCoarseIfFinite}) $ad(X)=\dim^{coa}_{rse}(X)=m$.
Suppose $x\in X$. Without loss of generality we may assume $x\in B$.
There is $s\in S$ such that $B(x,2M)\cap B\subset W_s$.
Hence $B(x,M)\subset B(W_s,M)\subset U_s$ and, since
any $y\in B(x,M)$ satisfies $dist(y,W_s)\leq d(y,x) < M < dist(y,B\setminus W_s)$,
we get $y\in e(W_s)$ which completes the proof.
\end{pf}

\section{Coarse dimension and absolute extensors}

In \cite{Dr$_1$} (Remark 2 on p.1097) Dranishnikov pointed out that $R_+$ is not an absolute
extensor in the category of proper metric spaces and coarse functions. He characterized proper metric spaces
of coarse dimension at most $n$ as those for which $R^{n+1}$
is an absolute extensor in the category of proper approximately
Lipschitz functions (Definition 4 on p.1105 and Theorem 6.6 on p.1111). 
That still left the door open to the possibility
of characterizing coarse dimension via $R^{n+1}$ being an absolute
extensor in the proper coarse category. The following result
clarifies that issue in negative.

\begin{Thm} \label{AbsExtensorsAndCoarseDimZero}
For a metric space $X$ the following conditions
are equivalent:

1. The coarse dimension of $X$ is at most $0$.

2. $Y$ is an absolute extensor of $X$ in the proper coarse category
for all $Y$.

3. $R_+$ is an absolute extensor of $X$ in the proper coarse category.
\end{Thm} 
\begin{pf} 1$\implies$2. It suffices to show that any unbounded
subset $A$ of $X$ is a coarsely proper and coarse retract of $X$.
Pick $x_0\in X$.
 Define by induction on $n$ an increasing sequence
$M_n$ of natural numbers and covers ${\cal {U}}^n$ of $X$ satisfying
the following properties:

a. $M_1=1$.

b. ${\cal {U}}^n$ is $M_n$-disjoint, the diameters of its elements
are smaller than $M_{n+1}$, and $L({\cal {U}}^n,X) > M_n$.

For each $U\in {\cal {U}}^n$ so that $U\cap A\ne\emptyset$,
pick $x_U\in U\cap A$ satisfying 
$d_X(x_U,x_0) > \sup\{d_X(x,x_0) \mid  x\in U\cap A\} -1/n$.

By induction on $n$ define a sequence of subsets $A_n$ of $X$
and a sequence of functions $r_n:A_n\to A$ as follows:

i. $A_1=A$ and $r_1=id_A$.

ii. $A_{n+1}$ is the union of those elements of ${\cal {U}}^{n+1}$
that intersect $A$.

iii. If $x\in U\setminus A_n$ and $U\cap A\ne\emptyset$
for some $U\in{\cal {U}}^{n+1}$, then $r_{n+1}(x)=x_U$.

Notice $X=\bigcup\limits_{n=1}^\infty A_n$ and let $r:X\to A$
be obtained by pasting all $r_n$. Observe that
$x\in U\in{\cal {U}}^k$ and $U\cap A\ne\emptyset$ implies $r(x)\in U$. Indeed, for each
$n$ there is a unique element $U_x^n\in {\cal {U}}^n$ containing $x$
and $U^i_x\subset U^j_x$ if $i < j$. Find the smallest
number $m$ so that $x\in A_m$. In that case $r(x)\in U^m_x$
by definition and $k$ must be at least $m$ so $U^m_x\subset U^k_x=U$.

We will show that $r$ is coarse by proving
$d_X(x,y) < M_n$ implies $d_X(r(x),r(y) \leq M_{n+2}$. 
 Indeed, if $d_X(x,y) < M_n$, then 
one of the following cases occurs:

Case 1. $U\cap A_n=\emptyset$, where $U$ is the unique element of ${\cal {U}}^{n+1}$
containing both $x_n$ and $y_n$.

Case 2. $U\cap A_n\ne\emptyset$, where $U$ is the unique element of ${\cal {U}}^{n+1}$
containing both $x_n$ and $y_n$.

In Case 1 the values $r(x)$ and $r(y)$ are identical.
In Case 2 both $r(x)$ and $r(y)$ belong to $U\cap A$ and
the set $U\cap A$ is of diameter at most $M_{n+2}$, so
$d_X(r(x),r(y) \leq M_{n+2}$ holds.

If $r$ is not coarsely proper, then there is a sequence $x_n\to\infty$ such that
$r(x_n)$ is bounded. Obviously, $x_n\notin A$ for almost all $n$.
Consider an element $U\in {\cal {U}}^k$ containing all of $r(x_n)$.
The way functions $r_m$ were defined implies that there is a sequence
of elements $U_n\in {\cal {U}}^{\alpha (n)}$ with $\alpha (n)\to\infty$
and all $U_n$ containing $U$, such that $U_n\cap A$ is of almost the same
diameter as $U\cap A$. That contradicts $A$ being unbounded.

2$\implies$3 is obvious.

3$\implies$1. Suppose $\dim^{coa}_{rse}(X) > 0$. 
By \ref{CharOfZeroAsympDim} there exists a number $M>0$ and a coarsely proper sequence 
$\{(x_n,y_n)\}_{n=1}^\infty$ of pairs of points in $X$ such that 
$dist(x_n,y_n)\to\infty$ and the points $x_n$ and $y_n$ can be 
$M$-scale connected in $X\setminus B(x_0,n)$
by long chain of length $L_n$  so that $L_n\to\infty$ as $n\to\infty$. 
We may assume $d_X(x_{n+j},x_n) >n$ and $d_X(y_{n+j},y_n) >n$
for all $n,j\ge 1$.
Let $B=\{x_n\}\cup\{y_n\}$. Define $f:B\to R_+$ by sending
$x_n$ to $n$ and $y_n$ to $n+n\cdot L_n$. Notice $f$ is coarsely proper and coarse.
Suppose $f$ extends to a coarse function $g:X\to R_+$.
Find $K > 0$ such that $d_X(x,y) \leq M$ implies $d(f(x),f(y))\leq K$.
Since $x_n$ and $y_n$ can be connected by a chain of $L_n$ points,
with consecutive points being separated by at most $M$,
$L_n\cdot n+n-n=d(f(x_n),f(y_n))\leq L_n\cdot K$ which leads to a contradiction
for $n > K$.
\end{pf}

\section{Open problems}

In \cite{Dr$_1$} (Problem 1 on p.1126) it is asked if
the asymptotic dimension of a proper metric space $X$ equals
the covering dimension of its Higson corona. Here is our version of that problem.

\begin{Problem}\label{ProblemAsdimInfCoarseFinite}
Is there a metric space $X$ of infinite asymptotic dimension and finite coarse dimension?
\end{Problem}

\begin{Problem}\label{ProblemMajorCoarseInfCoarseFinite}
Is there a metric space $X$ of infinite major coarse dimension and finite coarse dimension?
\end{Problem}

\begin{Def}\label{BoundedGeometryDef}
A metric space $X$ is of {\it bounded geometry}
if for every $M > 0$ there is a uniformly bounded cover ${\cal {U}}$
of $X$ of finite multiplicity and the Lebesque number at least $M$.
\end{Def}

\begin{Def}[\cite{Dr$_1$},p.1005]\label{Functiond(M)Def}
Suppose $X$ is a metric space of bounded geometry.
Given $M > 0$ let $d(M)=m({\cal {U}})-1$, where ${\cal {U}}$
 is a uniformly bounded cover ${\cal {U}}$ of minimal multiplicity
among those of the Lebesque number at least $M$.
$X$ is of {\it slow dimension growth} if $\lim\limits_{M\to\infty}\frac{d(M)}{M}=0$.
\end{Def}
Just as in \cite{Dr$_1$} (Problem 6 on p.1126) one can ask variants of
problems \ref{ProblemAsdimInfCoarseFinite} and \ref{ProblemMajorCoarseInfCoarseFinite} 
for spaces of bounded geometry or slow dimension growth.

\begin{Problem}\label{ProblemAsdimInfCoarseFiniteSlowGrowth}
Suppose $X$ is of slow dimension growth and finite coarse dimension.
Is asymptotic dimension of $X$ finite?
\end{Problem}

\begin{Problem}\label{ProblemMajorCoarseInfCoarseFiniteSlowGrowth}
Suppose $X$ is of slow dimension growth and finite coarse dimension.
Is the major coarse dimension of $X$ finite?
\end{Problem}

The above problems remain open for minor asymptotic dimension. All of the above
problems are of interest in case of $X$ being a finitely generated group
with word metric, especially $CAT(0)$ groups.

\begin{Problem}\label{ProblemDimOfProduct}
It is stated in \cite{DJ} that
$asdim(X \times Y) \leq asdim(X) + asdim(Y)$.
Are the corresponding results true for other dimensions?
\end{Problem}

\end{document}